\documentclass[12pt]{article}
\usepackage[utf8]{inputenc}
\usepackage[margin=1.25in]{geometry}
\usepackage{amsthm}
\usepackage{amsfonts}
\usepackage{amssymb}
\usepackage{amsthm}
\usepackage{hyperref}
\usepackage{tikz}
\usepackage{float}
\usepackage{microtype}

\newtheorem{theorem}{Theorem}[section]

\newtheorem{lemma}[theorem]{Lemma}
\newtheorem{definition}[theorem]{Definition}
\newtheorem{question}[theorem]{Question}
\newtheorem{proposition}[theorem]{Proposition}

\title{Logical limit laws for layered permutations and related structures}
\author{Samuel Braunfeld \and Matthew Kukla}
\date{}
\begin{document}

\begin{center}
	{\large \bf Logical limit laws for layered permutations and related structures}\\[10pt]
	
	Samuel Braunfeld$^\ast$ and Matthew Kukla$^\dag$\\[20pt]

	\footnotesize {\it $^\ast$Department of Mathematics, University of Maryland, College Park\\
		4176 Campus Dr, College Park, MD, USA 20742 \\
		Email: sbraunf@umd.edu}\\[10pt]
	
	\footnotesize {\it $^\dag$Department of Mathematics, University of Maryland, College Park\\
		4176 Campus Dr, College Park, MD, USA 20742 \\
		Email: mkukla1@umd.edu}\\[20pt]
\end{center}

\setcounter{page}{1} \thispagestyle{empty}

\baselineskip=0.30in

\normalsize

\abstract{We show that several classes of ordered structures (namely, convex
linear orders, layered permutations, and compositions) admit
first-order logical limit laws.
\bigskip


}

\section{Introduction}

Given a class of arbitrarily large finite structures, it is a classic
problem to ask,``What does a large random structure from this class look like?" Zero-one laws are one method of
gaining insight into this question for the case of properties expressible in
first-order logic.  A class $\mathcal{C}$ of structures is said to admit a zero-one law if, given a sentence $\varphi$ of first-order logic, 
the probability that a randomly selected $\mathcal{C}$-structure
of size $n$ satisfies $\varphi$ converges asymptotically to zero or one as $n$ goes 
to infinity.  For example, in the class of finite graphs (in the language with the edge relation), properties such as containing a specified subgraph or induced
subgraph are expressible in first-order logic, whereas properties such as
connectedness or Hamiltonicity are not.  It is a seminal result that this class
admits a zero-one law (\cite{fagin1976probabilities},\cite{glebskii1969range}).

The requirement that such a probability for a class of structures converges to
either zero or one is rather strict, and in general, not many classes admit a 
zero-one law. This is particularly evident when dealing with ordered structures.  
For example, consider the class of finite linear orders with two colors: 
red and blue.  The probability that a randomly selected structure's first point is 
red (or blue) is $1/2$.  A class of structures is said to admit a {\em logical 
limit law} when the probability that a randomly selected structure of size $n$ 
satisfies any given first order property converges as $n \to \infty$, so a zero-one law is the special case where the limiting probability is always zero or one.

We prove logical limit laws for various classes of ordered
structures.  In particular, we show the following result.

\begin{theorem}\label{theorem:intro}
Convex linear orders and layered permutations admit both unlabeled and labeled
limit laws.  Compositions admit an unlabeled limit law.
\end{theorem}

We first prove a logical limit law for convex linear orders using an adaptation of
Ehrenfeucht's proof for colored linear orders, which itself uses stationary states of
Markov chains.  We then use some basic definability results to transfer the limit law to layered permutations and compositions.

We now discuss some connections to previous work. First-order properties of finite permutations (when viewed as pairs of linear 
orders) were studied in \cite{Albert}. There the existence of a zero-one law was
disproven, and it was asked whether or not permutations could admit a logical
limit law; the answer to this turns out to be negative as well (as shown in \cite{Foy}). However, there is a rich study of various permutation classes other than the
class of all permutations, and we can also ask whether these admit logical limit
laws, as has been done for various restricted graph classes (e.g
\cite{kopczynski2019logical}).   Layered permutations are a
simple and well-studied example of such a permutation class.

In \cite{Cameron}, the limiting probability distributions of several
properties are computed for random compositions (what we call ``compositions,'' \cite{Cameron} refers to as ``preorders'').  Our result complements this by showing that the limiting probability of any first-order
property converges, although we do not describe any sort of limiting distribution.

\section{Convex Linear Orders}
\begin{definition}
The language of convex linear orders. $\mathcal{L}$, consists of two 
binary relations: $<$ (a total order on points) and $E$ (an equivalence relation). 

A {\em convex linear order} is an $\mathcal{L}$-structure $\mathfrak{C}$ such that $E$-classes are $<$-intervals (i.e., for $x, y, z \in \mathfrak{C}$ with  $x\mathrel{E} y$ and $x < z<y$, it holds that $z \mathrel{E} x,y$).
\end{definition}

In this section, we prove a logical limit law for the class of all finite convex
linear orders, closely following Ehrenfeucht's argument for colored linear orders as presented in Chapter 10 of  \cite{RG}.

We denote the convex linear order with  one point by $\bullet$.
\subsection{Preliminaries}
\begin{definition}\label{def:convlinorder}
Let $\mathfrak{C}$ be a convex linear order.  Define 
$\widehat{\mathfrak{C}}$ to be the convex linear order obtained by
adding one additional point to the last class of $\mathfrak{C}$.
\end{definition}

\begin{definition}\label{def:oplus}
For convex linear orders $\mathfrak{C}, \mathfrak{D}$, define $\mathfrak{C}
\oplus \mathfrak{D}$ as the convex linear order placing
$\mathfrak{D}$ $<$-after  $\mathfrak{C}$.
\end{definition}

Clearly, the last class of any convex linear order either contains exactly one
element or more than one element. This leads naturally to:

\begin{lemma} \label{lemma:con}
Every finite convex linear order containing $n$ points can be uniquely
constructed by applying $\widehat{(-)}$ and $- \oplus \bullet$ to $\bullet$ 
repeatedly.  This construction is done in $n-1$ steps.
\end{lemma}

\begin{proof}
We proceed by induction.  Let $\mathfrak{C}$ be a convex linear order of size $n$. If $n=1$, $\mathfrak{C} \simeq \bullet$.
  
Now assume that for arbitrary $n$, any convex linear order of size $n-1$ can be 
uniquely constructed from the operations above, in $n-2$ steps.  Let $\mathfrak{B}$ be $\mathfrak{C}$ with the $<$-last point removed.  Either the last class of 
$\mathfrak{C}$ contains one point, or it contains more than one point.  
If the former is true, then $\mathfrak{C} \simeq \mathfrak{B} \oplus \bullet$. 
In this case, 
$\mathfrak{C}$ cannot be obtained from $\widehat{(-)}$ from a structure of size $n-1$, as the last class of $\mathfrak{C}$ needs 
to contain exactly one point. If the latter is true, then $\mathfrak{C} \simeq 
\widehat{\mathfrak{B}}$, and cannot be obtained from applying $- \oplus \bullet$ to a structure of size $n-1$.
\end{proof}

From Lemma \ref{lemma:con}, we see that starting with $\bullet$ and randomly applying $- \oplus \bullet$
and $\widehat{(-)}$ (each with probability $1/2$) $n-1$ times  will uniformly randomly sample all
possible convex linear orders of size $n$. 

\begin{definition}\label{def:ef}
Suppose $\mathfrak{A}$ and $\mathfrak{B}$ are structures in a
language consisting of binary relations $R_1, \ldots, R_n$.  An {\em
Ehrenfeucht–Fraïssé game} (sometimes referred to as a {\em ``back and forth game''}) of length $k$ between $\mathfrak{A}$ and $\mathfrak{B}$  
is a game between two players (referred to as Duplicator and
Spoiler).  In each round, Spoiler plays by selecting and marking a point on
either structure.  Duplicator responds by marking a corresponding point on the
structure which Spoiler did not choose from.  After $k$ rounds, the points 
$x_1, \ldots, x_k$ have been marked on $\mathfrak{A}$, and $y_1, \ldots, y_k$ on
$\mathfrak{B}$.  Duplicator has won if the map  sending each $x_i$ to the 
corresponding $y_i$ is an isomorphism (i.e., $x_1, \ldots, x_k, y_1, \dots, y_k$ satisfy $x_i 
\mathrel{R}_\ell x_j \iff y_i \mathrel{R}_\ell y_j$ for all $i,j \in [k]$ and $\ell \in [n]$). Otherwise, Spoiler has won. 
\end{definition}

The quantifier depth of a first order sentence counts the maximum depth of
nested quantifiers.  A formal definition via induction is found on page 16
of \cite{RG}.

\begin{theorem}
In the Ehrenfeucht–Fraïssé game of length $k$ between $\mathfrak{A}$ and $\mathfrak{B}$, Duplicator has a winning strategy iff $\mathfrak{A}$ and $\mathfrak{B}$ agree
on all sentences of quantifier depth at most $k$.
\end{theorem}

\begin{proof}
A proof can be found in Section 2.3.1 of \cite{RG} for the case of graphs; the
general case follows an essentially identical argument.
\end{proof}

We write $\mathfrak{A} \equiv_k \mathfrak{B}$ when $\mathfrak{A}$ and $\mathfrak{B}$ agree on all sentences of quantifier depth at most $k$.

\begin{lemma} \label{lemma:add}
Let $\mathfrak{M}, \mathfrak{N}, \mathfrak{M}', \mathfrak{N}'$ be convex linear orders such that $\mathfrak{M} \equiv_k \mathfrak{N}$ and 
$\mathfrak{M}' \equiv_k \mathfrak{N}'$.  Then, $\mathfrak{M} \oplus
\mathfrak{M}' \equiv_k \mathfrak{N} \oplus \mathfrak{N}'$.
\end{lemma}

\begin{proof}
We will show that in any Ehrenfeucht–Fraïssé game of length $k$, Duplicator has a 
winning strategy. Consider such a game between $\mathfrak{M} \oplus \mathfrak{M}'$ 
and $\mathfrak{N} \oplus \mathfrak{N}'$. 
If Spoiler chooses any element in $\mathfrak{M}$ (respectively, $\mathfrak{M}'$), 
then Duplicator responds as they would in a length $k$ game between
$\mathfrak{M}$ and $\mathfrak{N}$ (respectively, between $\mathfrak{M}'$ and
$\mathfrak{N}'$), and vice-versa if Spoiler chooses a point in $\mathfrak{N}
\oplus \mathfrak{N}'$. We will show this gives Duplicator a winning strategy in
the Ehrenfeucht–Fraïssé game between $\mathfrak{M} \oplus \mathfrak{M}'$ and
$\mathfrak{N} \oplus \mathfrak{N}'$.  
Let $A \subset \mathfrak{M}, A' \subset \mathfrak{M}', B \subset \mathfrak{N}, B' \subset \mathfrak{N}'$ be the elements chosen in the game. Then $A \simeq B$ and $A' \simeq B'$. But it is easy to see that $\oplus$ preserves isomorphism, so $A \oplus A' \simeq B \oplus B'$.
\end{proof}

A very similar argument gives the following:

\begin{lemma} \label{lemma:hat}
Suppose $\mathfrak{M} \equiv_k \mathfrak{N}$, then, $\widehat{\mathfrak{M}} \equiv_k 
\widehat{\mathfrak{N}}$.
\end{lemma}

\begin{proof}
We again show that in an Ehrenfeucht–Fraïssé game of length $k$, Duplicator has a 
winning strategy.  For any move by Spoiler in $\mathfrak{M}$ (or $\mathfrak{N}$),
Duplicator responds by playing as they would normally would in an
Ehrenfeucht–Fraïssé game between $\mathfrak{M}$ and $\mathfrak{N}$; because 
$\mathfrak{M} \equiv_k \mathfrak{N}$, Duplicator always has a winning move in 
response to Spoiler for any such play.  If Spoiler plays the last point in the 
last class of $\widehat{\mathfrak{M}}$ or $\widehat{\mathfrak{N}}$ (that is, the 
point added by $\widehat{(-)}$), Duplicator can always respond with the 
corresponding point at the end of $\widehat{\mathfrak{N}}$ or $\widehat{\mathfrak{M}}$.  Hence, Duplicator has a response for any of Spoiler's moves in a length $k$
Ehrenfeucht–Fraïssé game between $\widehat{\mathfrak{M}}$ and $\widehat{\mathfrak{N}
}$, so $\widehat{\mathfrak{M}} \equiv_k \widehat{\mathfrak{N}}$.
\end{proof}

\begin{lemma} \label{lemma:2k}
For two finite linear orders $N, M$ having $n$ and $m$ points respectively, 
$N \equiv_k M$ iff $n = m$ or $n, m \geq 2^k-1$.
\end{lemma}

\begin{proof}
This is Lemma $2.6.3$ in \cite{RG}.
\end{proof}

\begin{lemma} \label{lemma:1024}
For a convex linear order $\mathfrak{M}$ and $k \in \mathbb{N}$, there exists $\ell \in
\mathbb{N}$ such that for all $s, t > \ell$, \[ \bigoplus_s \mathfrak{M} \equiv_k
\bigoplus_t \mathfrak{M}\]
\end{lemma}

\begin{proof}
We reduce this to a case of the previous lemma. Let $\ell = 2^{k-1}$.  Let $O_s$ be a linear order 
with $s$ points, each corresponding to a copy of $\mathfrak{M}$ in
$\bigoplus_s \mathfrak{M}$, and define $O_t$ likewise for $t$ and
$\bigoplus_t \mathfrak{M}$.  In an  Ehrenfeucht–Fraïssé game of length $k$ between
$\bigoplus_s \mathfrak{M}$ and $\bigoplus_t \mathfrak{M}$, we will show that
Duplicator has a winning strategy.  If Spoiler picks a point in the $i$th copy
of $\mathfrak{M}$ in $\bigoplus_s \mathfrak{M}$, we view this as 
Spoiler picking the $i$th point in $O_s$ if it were playing a length-$k$
Ehrenfeucht–Fraïssé game between $O_s$ and $O_t$.  By Lemma \ref{lemma:2k}, Duplicator has a 
response in $O_t$; suppose this response is the $j$th point.  To have a winning 
strategy in the Ehrenfeucht–Fraïssé game between $\bigoplus_s \mathfrak{M}$ and 
$\bigoplus_t \mathfrak{M}$, Duplicator can select the same point in $\mathfrak{M}$ 
which Spoiler selected, but in the $j$th copy of $\mathfrak{M}$ in $\bigoplus_t 
\mathfrak{M}$.
\end{proof}

\subsection{The limit law}
It is in general important to note the distinction between labeled
and unlabeled limit laws, which count structures of size $n$ differently.  
Labeled limit laws count all possible structures
over the universe $\{1, \ldots, n\}$ as $n \to \infty$, whereas unlabeled limit laws count all
structures up to isomorphism over $\{1, \ldots, n\}$. In general, a labeled
limit law for a given class does not imply an unlabeled limit law for 
the class, and vice-versa.  However, as finite linearly ordered structures have
no nontrivial automorphisms, labeled and unlabeled limit laws are equivalent. Thus we will not distinguish between them in the section or the next.

Given a first-order sentence $\varphi$ having quantifier rank $k$, we compute
the limiting probability of $\varphi$ by associating to it a Markov chain
$M_\varphi$.  For a $\equiv_k$-class of convex linear orders $C$, and for some/any $\mathfrak{M} \in C$, we make the following definitions.\[ C \oplus \bullet := \left [ \mathfrak{M}
\oplus \bullet \right ]_{\equiv_k} \]

and \[\widehat{C} := \left [ \widehat{\mathfrak{M}} \right ]_{\equiv_k} \]

These operations are well-defined, as any choice of $\mathfrak{M}$ yields a $\equiv_k$-equivalent result by Lemmas \ref{lemma:add} and \ref{lemma:hat}.

The states of $M_\varphi$ are $\equiv_k$-classes
of $\mathcal{L}$-structures (where $k$ is the quantifier depth of $\varphi$); there are finitely many of such classes, by
Theorem 2.2.1 of \cite{RG}.  For a $\equiv_k$-class $C$, there are two possible transitions out of $C$: one to 
$C \oplus \bullet$, and one to $\widehat{C}$, each having probability $1/2$.  
The starting state of $M_\varphi$ is $\bullet$ (we slightly abuse notation by writing $\bullet$ to also 
mean $\left [ \bullet \right ]_{\equiv_k}$). 

\begin{definition}\label{def:periodic} 
A Markov chain $M$ is {\em fully aperiodic} if there do not exist 
disjoint sets of $M$-states $P_0, P_1, \ldots, P_{d-1}$ for some $d > 1$ such 
that for every state in $P_i$, $M_\varphi$ transitions to a state in $P_{i+1}$ 
with probability $1$ (with $P_{d-1}$ transitioning to $P_0$), i.e every
state is aperiodic.
\end{definition}

We next state a variant of the fundamental theorem of Markov chains which does not
assume irreducibility.

\begin{proposition} \label{prop:markov}  Let $M$ be a finite, fully aperiodic 
Markov chain with initial state $S$, and let $Pr^{n-1}(S, Q)$ denote the 
probability that $M$ is in state $Q$ after $n-1$ steps.  Then, for any
$Q$, $\lim_{n \to \infty} Pr^{n-1}(S, Q)$ converges. 
\end{proposition}

\begin{proof}
This is contained in the discussion following Theorem 0.3.1 of \cite{Markov}.
\end{proof}

\begin{lemma} \label{lemma:fully-aperiodic}
For any first-order sentence $\varphi$, $M_\varphi$ is fully aperiodic.
\end{lemma}

\begin{proof}
Suppose $M_\varphi$ were not fully aperiodic. Then, there would exist disjoint sets of 
$M$-states  ($\equiv_k$-classes) $P_0, P_1, \ldots, P_{d-1}$ for some $d > 1$ such that for every
state in $P_i$, $M_\varphi$ transitions to a state in $P_{i+1}$ with probability
$1$ (with $P_{d-1}$ transitioning to $P_0$).  Write $i \bullet$ to mean
$\bigoplus_i \bullet$.  Thus, for any $Q \in P_0$, $Q \oplus
i \bullet$ is in $P_0$ iff $d \mid i$.  But by
Lemma \ref{lemma:add} and Lemma \ref{lemma:1024}, $Q \oplus i \bullet
\equiv_k Q \oplus (i+1) \bullet$ for sufficiently large $i$, contradicting 
this claim. 
\end{proof}

\begin{theorem} \label{theorem:convlinorder}
Convex linear orders admit a logical limit law.
\end{theorem}

\begin{proof}
Fix a first-order sentence $\varphi$, and consider the Markov chain $M_\varphi$.
For each state $S$ in $M_\varphi$, either 
each structure in $S$ satisfies
$\varphi$ or no structures in  $S$ satisfy $\varphi$.  Let $S_\varphi$ denote the set of
states in $M_\varphi$ for which all structures in that
state satisfy  $\varphi$.  By the fact that $\widehat{(-)}$ and $- \oplus
\bullet$ are well-defined on $\equiv_k$-classes (it does not matter which
structure in the class is chosen), we can view moving $n-1$ steps in $M_\varphi$ as
starting with any structure in the current state, applying $\widehat{(-)}$ or 
$- \oplus \bullet $ to it $n-1$ times,  and taking the resulting $\equiv_k$-class at the end. Thus by the comment after Lemma \ref{lemma:con}, moving $n-1$ steps from the starting state $\bullet$ is the same as uniformly randomly picking a size $n$ 
structure and then taking its $\equiv_k$-class.
Therefore, the probability that after $n$ steps, the chain is in a state of 
$S_\varphi$ is same as probability that uniformly randomly selected structure of
size $n$ satisfies $\varphi$.
\if{false}
By the comment after Lemma \ref{lemma:con}, the probability that a uniformly, 
randomly selected structure of size $n$ satisfies $\varphi$ is precisely the 
probability that, after $n-1$ steps, $M_\varphi$ is in a state in $S_\varphi$.  
This probability is given by
\[ \sum_{Q \in S_\varphi}^{}Pr^{n-1}(\bullet, Q) \]
where $Pr^{n-1}(\bullet, Q)$ denotes the transition probability from
$\bullet$ to $P$ after $n$ steps . 
\fi
So it suffices to show that $\lim_{n \to \infty} \sum_{Q \in
S_\varphi}^{}Pr^{n-1}(\bullet, Q)$ converges. Because $M_\varphi$ has 
finitely many states, \[\lim_{n \to \infty}  \sum_{Q \in
S_\varphi}^{}Pr^{n-1}(\bullet, Q) = \sum_{Q \in
S_\varphi}^{}\lim_{n \to \infty}Pr^{n-1}(\bullet, Q) \]

It now suffices to show that $\lim_{n \to \infty}
Pr^n(\bullet, Q)$ exists for every state $Q$ of $M_\varphi$.  But this follows from Proposition \ref{prop:markov}.
\end{proof}

\section{Layered Permutations}
Permutations can be viewed as structures in a language with two linear orders,
$<_1$, $<_2$. The order $<_1$  gives the unpermuted order of the points (before 
applying the permutation) and $<_2$ describes the points in permuted order.  An embedding of one such structure into another then corresponds to the usual notion of pattern containment.

The properties of permutations expressible by a first-order sentence in this language are explored in Section 3 of \cite{Albert}. These include the containment and avoidance of (generalized) patterns, concepts related to the substitution decomposition, and sortability properties such as $k$-stack sortability.

\begin{definition}\label{def:lperm}
	Given a
	permutation $P$, a {\em block} is a maximal subset $B \subset P$ that is an interval with respect to $<_1$ and $<_2$, and is monotone.
	
A {\em layered permutation} consists of increasing blocks of decreasing permutations.
\end{definition}
\begin{figure}[H] \label{figure:lperm}
    \centering
    \caption{Illustration of a layered permutation.  From left to right, $<_1$
    is increasing; from bottom to top, $<_2$ is increasing.}
    \label{fig:label1}
\begin{tikzpicture}[x=0.75pt,y=0.75pt,yscale=-1,xscale=1]

\draw   (64,145) -- (114,145) -- (114,195) -- (64,195) -- cycle ;
\draw   (114,95) -- (164,95) -- (164,145) -- (114,145) -- cycle ;
\draw   (164,45) -- (214,45) -- (214,95) -- (164,95) -- cycle ;
\draw    (36,220) -- (249,220) ;
\draw [shift={(251,220)}, rotate = 180] [color={rgb, 255:red, 0; green, 0; blue, 0 }  ][line width=0.75]    (10.93,-3.29) .. controls (6.95,-1.4) and (3.31,-0.3) .. (0,0) .. controls (3.31,0.3) and (6.95,1.4) .. (10.93,3.29)   ;
\draw    (36,220) -- (36,19) ;
\draw [shift={(36,19)}, rotate = 450] [color={rgb, 255:red, 0; green, 0; blue, 0 }  ][line width=0.75]    (10.93,-3.29) .. controls (6.95,-1.4) and (3.31,-0.3) .. (0,0) .. controls (3.31,0.3) and (6.95,1.4) .. (10.93,3.29)   ;

\draw (66,150.4) node [anchor=north west][inner sep=0.75pt]    {$\bullet $};
\draw (97,175.4) node [anchor=north west][inner sep=0.75pt]    {$\bullet $};
\draw (116,100.4) node [anchor=north west][inner sep=0.75pt]    {$\bullet $};
\draw (150,128.4) node [anchor=north west][inner sep=0.75pt]    {$\bullet $};
\draw (132,113.4) node [anchor=north west][inner sep=0.75pt]    {$\bullet $};
\draw (183,60.4) node [anchor=north west][inner sep=0.75pt]    {$\bullet $};
\draw (128,225.4) node [anchor=north west][inner sep=0.75pt]    {$< _{1}$};
\draw (4,102.4) node [anchor=north west][inner sep=0.75pt]    {$< _{2}$};

\end{tikzpicture}
\end{figure}
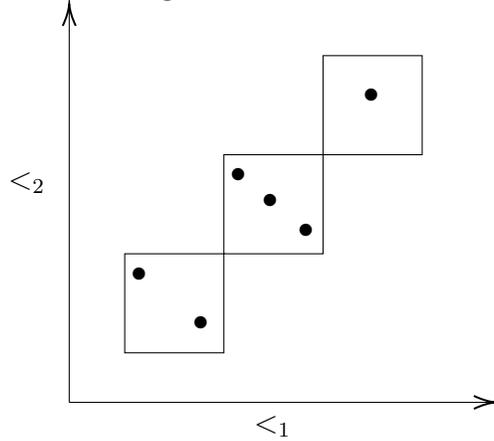

Note that in a layered permutation, two points are in the same block iff $<_1$ and $<_2$ disagree between them.

In many cases, it is useful to transfer a logical limit law on one class of
structures to another class of similar structures.  
This is possible when the following interdefinability condition is satisfied. (This is surely well-known, but we could not find an explicit statement.)

\begin{definition}(Uniform Interdefinability) \label{id}\\ 
Let $\mathcal{L}_0$, $\mathcal{L}_1$ be languages, and $\mathcal{C}_0$, 
$\mathcal{C}_1$ be classes of finite $\mathcal{L}_0$, $\mathcal{L}_1$ structures 
respectively, with a common domain of $\left [ n  \right ] = \left \{ 1,\cdots,n
\right \}$.  $\mathcal{C}_0$ and $\mathcal{C}_1$ are said to be 
{\em uniformly interdefinable} if there exists a map $f_I: \mathcal{C}_0 \to 
\mathcal{C}_1$ (which is a bijection on structures of size $n$, for all $n$), along 
with formulae $\varphi_{R_{0,i}}, \varphi_{R_{1,i}}$ for each relation $R_{0,
i}$ in $\mathcal{L}_0$ 
and $R_{1, i}$ in $\mathcal{L}_1$ respectively such that, 
for each $\mathfrak{M}_0$ in $\mathcal{C}_0$ and $\mathfrak{M}_1$ in $\mathcal{C}_1$:

\begin{itemize}
    \item $\mathfrak{M}_0 \models R_{0,i}(\bar{x}) \iff 
    f_I(\mathfrak{M}_0) \models \varphi_{R_{0,i}}(\bar{x})$
    
    \item $\mathfrak{M}_1 \models R_{1,i}(\bar{x}) \iff f_I^{-1}(\mathfrak{M}_1) 
    \models \varphi_{R_{1,i}}(\bar{x})$
\end{itemize}

\end{definition}
Although uniform interdefinability may appear to be asymmetric, it is in fact a
symmetric relation, using that $f_I$ is bijective.

\begin{lemma} \label{lemma:translation}
Let $\mathcal{L}_0$, $\mathcal{L}_1$ be languages, $\mathcal{C}_0$,
$\mathcal{C}_1$ classes of $\mathcal{L}_0$, $\mathcal{L}_1$-structures
respectively, $f$ a map from the set of  $\mathcal{L}_0$-structures to the set
of $\mathcal{L}_1$-structures, and $g$ a map from the set of 
$\mathcal{L}_0$-sentences to the set of $\mathcal{L}_1$-sentences such that, for
any $\mathcal{C}_0$-structure $\mathfrak{M}$ and $\mathcal{L}_0$-sentence
$\varphi$:

\begin{enumerate}
    \item $\mathfrak{M} \models \varphi \iff 
    f(\mathfrak{M}) \models g(\varphi)$ \item $f$ is a bijection between
    $\mathcal{C}_0$ and $\mathcal{C}_1$ structures of size $n$
    \item The class $\mathcal{C}_1$ admits a logical limit law
\end{enumerate}

Then, $\mathcal{C}_0$ admits a logical limit law as well.
\end{lemma}

\begin{proof}
Let $\varphi$ be a sentence in $\mathcal{L}_0$ and $a_0$ the number of $\mathcal{C}_0$-structures of size $n$ satisfying $\varphi$.  
Likewise, let $a_1$ be the number of $\mathcal{L}_1$-structures of size $n$ which satisfy $g(\varphi)$.
For a randomly selected $\mathcal{C}_0$-structure $\mathfrak{M}$ (of size $n$),
the probability that $\mathfrak{M} \models \varphi$ is $\frac{a_0}{\left | 
\mathcal{C}_0 \right |}$, and the probability that $f(\mathfrak{M}) \models g(\varphi)$ in $\mathcal{C}_1$ is $\frac{a_1}{\left | \mathcal{C}_1 \right |}$.  Since $f$ is a
bijection on structures of size $n$, we have that $\left | \mathcal{C}_1 \right | = \left | \mathcal{C}_0 \right |$ and together with (1), that $a_1 = a_0$.  Hence, the probabilities are equal for any 
$\varphi$; because $\mathcal{C}_1$ admits a limit law, 
$\mathcal{C}_0$ admits a limit law as well.
\end{proof}

\begin{lemma}\label{lemma:inter-translation}
Let $\mathcal{C}_0$, $\mathcal{C}_1$ be uniformly interdefinable classes of
$\mathcal{L}_0$, $\mathcal{L}_1$ structures.  If $\mathcal{C}_1$ admits a
logical limit law, $\mathcal{C}_0$ admits one as well.
\end{lemma}

\begin{proof}
We show that the maps $f$, $g$ exist as in Lemma \ref{lemma:translation}.  Take $f=f_I$ and $g$ as the map 
which sends an $\mathcal{L}_0$ sentence to the $\mathcal{L}_1$ sentence where each 
occurrence of $R_{0,i}$ is replaced by $\varphi_{R_{0,i}}$. Given a $\mathcal{C}_0$-structure $\mathfrak{M}_0$
and an atomic $\mathcal{L}_0$-
formula $\varphi(\bar{x})$, and a tuple $\bar{m} \subset \mathfrak{M}_0$,  we have $\mathfrak{M}_0 \models \varphi(\bar{m}) \iff f_I(\mathfrak{M}_0) \models g(\varphi)(\bar{m})$ by the definition 
of uniform interdefinability. When $\varphi$ is
nonatomic, the same statement follows from a standard induction on the complexity of
$\varphi$.  
\if{false}
$$\mathfrak{M}_0 \models R_0(\bar{x}) \iff f_I(\mathfrak{M}_0) \models \varphi_{R_0}(\bar{x}) = 
g(R_0(\bar{x}))$$
and identically for $R_1$, $\varphi_{R_1}$.
\fi
Furthermore, by Definition \ref{id}, $f_I$ is a bijection on structures of size
$n$, and therefore, $f_I, g$ are as desired.
\end{proof}

Let $\mathcal{L}_0 = \{ <_1, <_2 \}$ and $ \mathcal{L}_1 = \{ E, < \}$ denote the 
languages of layered permutations and convex linear orders respectively, and 
let $\mathcal{C}_0, \mathcal{C}_1$ classes of isomorphism types of 
$\mathcal{L}_0, \mathcal{L}_1$-structures respectively.  We define a map
from layered permutations to convex linear orders which sends blocks of a layered permutation to convex equivalence classes, and
points in each block of a layered permutation to points in the same equivalence
class such that $<_1$ agrees with $<$ (see Figure 2).  Formally, this is a map 
$f_I\colon \mathcal{C}_0 \to \mathcal{C}_1$ is defined such that for 
$\mathfrak{M}_0$ in $\mathcal{C}_0$ and $\mathfrak{M}_1$ in $\mathcal{C}_1$:
\begin{itemize}
    \item $f_I(\mathfrak{M}_0) \models a < b \iff \mathfrak{M}_0 \models a <_1 b $
    \item $f_I(\mathfrak{M}_0) \models a \mathrel{E} b \iff \mathfrak{M}_0 \models
    (a <_1 b \wedge a >_2 b) \vee (b <_1 a \wedge b >_2 a)$
\end{itemize}
The relations $<_1$, $<_2$ in the language of layered permutations are rewritten
in the language of convex equivalence relations using the following
rules on atomic formulas:
\begin{itemize}
    \item $\varphi_{<_1}:a <_1 b \rightsquigarrow a < b $
    \item $\varphi_{<_2}:a <_2 b \rightsquigarrow (a \mathrel{E} b \wedge b < a) 
    \vee (\neg (a \mathrel{E} b) \wedge a < b)$
\end{itemize}

\begin{figure}[H] \label{fig:2}


\tikzset{every picture/.style={line width=0.75pt}} 
\center

\begin{tikzpicture}[x=0.75pt,y=0.75pt,yscale=-1,xscale=1]
\draw   (299,115) -- (331,115) -- (331,145) -- (299,145) -- cycle ;
\draw   (332,84) -- (364,84) -- (364,114) -- (332,114) -- cycle ;
\draw   (301,28) .. controls (301,26.34) and (302.34,25) .. (304,25) -- (325,25) .. controls (326.66,25) and (328,26.34) .. (328,28) -- (328,37) .. controls (328,38.66) and (326.66,40) .. (325,40) -- (304,40) .. controls (302.34,40) and (301,38.66) .. (301,37) -- cycle ;
\draw   (337,28) .. controls (337,26.34) and (338.34,25) .. (340,25) -- (359,25) .. controls (360.66,25) and (362,26.34) .. (362,28) -- (362,37) .. controls (362,38.66) and (360.66,40) .. (359,40) -- (340,40) .. controls (338.34,40) and (337,38.66) .. (337,37) -- cycle ;
\draw    (306,122.5) -- (306,48.5) ;
\draw [shift={(306,46.5)}, rotate = 450] [color={rgb, 255:red, 0; green, 0; blue, 0 }  ][line width=0.75]    (10.93,-3.29) .. controls (6.95,-1.4) and (3.31,-0.3) .. (0,0) .. controls (3.31,0.3) and (6.95,1.4) .. (10.93,3.29)   ;
\draw    (321,132.5) -- (321,48.5) ;
\draw [shift={(321,46.5)}, rotate = 450] [color={rgb, 255:red, 0; green, 0; blue, 0 }  ][line width=0.75]    (10.93,-3.29) .. controls (6.95,-1.4) and (3.31,-0.3) .. (0,0) .. controls (3.31,0.3) and (6.95,1.4) .. (10.93,3.29)   ;
\draw    (348,97.5) -- (348,46.5) ;
\draw [shift={(348,44.5)}, rotate = 450] [color={rgb, 255:red, 0; green, 0; blue, 0 }  ][line width=0.75]    (10.93,-3.29) .. controls (6.95,-1.4) and (3.31,-0.3) .. (0,0) .. controls (3.31,0.3) and (6.95,1.4) .. (10.93,3.29)   ;

\draw (301,118.4) node [anchor=north west][inner sep=0.75pt]    {$\bullet $};
\draw (316,129.4) node [anchor=north west][inner sep=0.75pt]    {$\bullet $};
\draw (343,93.4) node [anchor=north west][inner sep=0.75pt]    {$\bullet $};
\draw (302,27.4) node [anchor=north west][inner sep=0.75pt]    {$\bullet $};
\draw (317,27.4) node [anchor=north west][inner sep=0.75pt]    {$\bullet $};
\draw (344,27.4) node [anchor=north west][inner sep=0.75pt]    {$\bullet $};
\draw (243,122.4) node [anchor=north west][inner sep=0.75pt]    {$ \begin{array}{l}
\mathfrak{M}_{0}\\
\end{array}$};
\draw (228,20.4) node [anchor=north west][inner sep=0.75pt]    {$f_I(\mathfrak{M}_{0})$};
\draw (289,72.4) node [anchor=north west][inner sep=0.75pt]    {$f_I$};

\end{tikzpicture}

\caption{Illustration of the map $f_I$.  Blocks of the layered permutation
$\mathfrak{M}_0$ are mapped to equivalence classes of $f_I(\mathfrak{M}_0)$, and
points are mapped in an order-preserving manner.}
\end{figure}
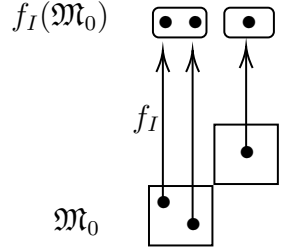       

It is perhaps clear from Figure 2 that convex linear orders and layered
permutations are uniformly interdefinable, but we now verify the details.

\begin{lemma} \label{lemma:id2}
Layered permutations and convex linear orders are uniformly interdefinable.
\end{lemma}

\begin{proof}
A finite layered permutation is determined, up to isomorphism, by the number of 
points in each of its blocks; likewise, a finite convex linear order is determined 
in the same manner by the number of points in each equivalence class.  Because
the map $f$ sends blocks of $\mathfrak{M}_0$ to equivalence classes of
$f(\mathfrak{M}_0)$ of the same size, $f$ is injective.  The number of possible
size $n$ convex linear orders is equal to the number of size $n$ layered
permutations (because every block structure is allowed), therefore, $f$ is
bijective as well.\\
Let $\mathfrak{M}_0$, $\mathfrak{M}_0'$ be two $\mathcal{C}_0$ -structures (layered
permutations) with $f(\mathfrak{M}_0) = f(\mathfrak{M}_0')$. From the definition of $<_1$, it is clear 
that $$\mathfrak{M}_0 \models a <_1 b \iff f(\mathfrak{M}_0) \models a < b$$
We verify $$\mathfrak{M}_0 \models a <_2 b \iff f(\mathfrak{M}_0) \models (a \mathrel{E} b \wedge b < a) \vee (\neg a \mathrel{E} b \wedge a < b)$$ \label{id:permutation} 

\noindent $(\Rightarrow)$
Suppose $a$ and $b$ are in the same block, then, because $\mathfrak{M}_0 \models
a <_2 b$, we have that $\mathfrak{M}_0
\models b <_1 a$, so $f(\mathfrak{M}_0) \models b < a$ (because the orders $<_1$
and $<$ agree).  When $a$ and $b$ are in
different blocks, $f(\mathfrak{M}_0) \models \neg a \mathrel{E} b$; furthermore,
because $\mathfrak{M}_0 \models a <_2 b$, and $a$, $b$ are in different blocks, 
$\mathfrak{M}_0 \models a <_1 b$, so $f(\mathfrak{M}_0) \models a < b$.\\

\noindent $(\Leftarrow)$
First suppose $a$ and $b$ are in the same equivalence class.
Then, $f(\mathfrak{M}_0) \models b < a$ and $\mathfrak{M}_0 \models b <_1 a$.
Because $a$, $b$ are in the same class, they are in the same block of
$\mathfrak{M}_0$. Since $\mathfrak{M}_0 \models b <_1 a$ and $a$, $b$ are
in the same block, we have $\mathfrak{M}_0 \models a <_2 b$.  

When $a$ and $b$ are in
different equivalence classes, $a < b$ in $\mathfrak{M}_0$, so $a <_1 b$ in $\mathfrak{M}_1$.  Since
$a$ and $b$ are in classes and thus different blocks, the orders $<_1$ and $<_2$ agree, giving $a <_2 b$.

\end{proof}

\begin{theorem}\label{theorem:permlimlaw}
Layered permutations admit a logical limit law.
\end{theorem}

\begin{proof}
By Lemma \ref{lemma:id2}, layered permutations are uniformly interdefinable with 
convex linear orders.  Because convex linear orders admit a logical limit law, 
layered permutations admit one as well by Lemma \ref{lemma:inter-translation}. 
\end{proof}

Zero-one laws have been extensively studied in the context of {\em homogeneous} structures, in the sense of \cite{cameron2002homogeneous}. The homogeneous permutations are classified in \cite{cameron2002homogeneous}, and Theorem \ref{theorem:permlimlaw} completes the picture of their logical limit laws. The increasing and decreasing permutations are uniformly interdefinable with linear orders, and so admit a zero-one law (e.g., see Section 2.6.2 of \cite{RG}), while the class of all permutations doesn't admit a logical limit law by \cite{Foy}. The remaining classes are the layered and ``skew layered'' permutations (consisting of decreasing blocks of increasing permutations), which are uniformly interdefinable by replacing $<_2$ with its reverse.

\section{Compositions}
Informally, a composition consists of an equivalence relation
$\mathrel{E}$ along with a linear order $\prec_1$ on $\mathrel{E}$-classes (but
not on points of the classes themselves).  This description naturally 
corresponds to a composition of $n$ in the usual sense, i.e. an ordered tuple of 
positive integers summing to $n$, and the embeddability order agrees with the order
on compositions in the usual sense given in \cite{bergeron1995standard}.

Compositions can be obtained from
convex linear orders by forgetting the order between points within the same 
$E$-class.  We formalize this by passing through the following notion of fractured 
orders.

Let $\mathcal{L}_0 = \{ \mathrel{E}, < \}$ be the language of convex linear orders 
as before.  Define a language $\mathcal{L}_1 = \{\mathrel{E}, \prec_1, \prec_2 \}$
consisting of three relation symbols (an equivalence relation and two partial
orders), and let the reduct $\mathcal{L}_2 \subset \mathcal{L}_1$ be given by 
$\mathcal{L}_2 = \{ \mathrel{E}, \prec_1 \}$.  In a fractured order, we start with a
convex linear order $<$ and break it into two parts: $\prec_1$ between
$\mathrel{E}$-classes, and $\prec_2$ within $\mathrel{E}$-classes.  
Formally, we define the class of finite fractured orders $\mathcal{F}$ 
to be the class of $\mathcal{L}_1$-structures satisfying:

\begin{enumerate}
	\item $\prec_1$, $\prec_2$ are partial orders
	\item $\mathrel{E}$ is an equivalence relation
    \item Distinct points $a$, $b$ are $\prec_1$-comparable iff they are not
    $E$-related
    \item Distinct points $a$, $b$ are $\prec_2$-comparable iff they are 
    $E$-related
    \item $a \mathrel{E} a'$, $a \prec_1 b \Rightarrow a' \prec_1 b$ (convexity)
\end{enumerate}

Although $\prec_1$ is a partial order on points, Axioms 3 and 5 say it is
essentially equivalent to a linear order on E-classes. So, compositions may be 
defined formally as  $\mathcal{L}_2$-reducts of 
fractured orders.  To show that compositions admit a logical limit law, we show 
that fractured orders are uniformly interdefinable with convex linear orders, 
and that any composition admits a unique expansion to a fractured 
order.

\begin{lemma} \label{lemma:c-oer}
	Convex linear orders and finite fractured orders are uniformly interdefinable.
\end{lemma}

\begin{proof}
	Define a map $f_I:\mathcal{F} \to \mathcal{C}_0$ such that:
	
	\begin{itemize}
		\item $\mathfrak{M}_1 \models a \mathrel{E} b \iff f_I(\mathfrak{M}_1)
		\models a \mathrel{E} b$ 
		
		\item $\mathfrak{M}_1 \models a \prec_1 b \iff f_I(\mathfrak{M}_1) \models 
		\neg a \mathrel{E} b \wedge a < b$
		
		\item $\mathfrak{M}_1 \models a \prec_2 b \iff f_I(\mathfrak{M}_1) \models a
		\mathrel{E} b \wedge a < b$
		
	\end{itemize}
	The order $<$ is total and convex by Axioms $3$ and $4$ in the definition of
    $\mathcal{F}$.  Both convex linear orders and compositions are determined 
	(up to isomorphism) by the number of points in each class, therefore the number of fractured orders of
	size $n$ is equal to the number of convex linear orders of
	size $n$.  Because $f_I$ preserves $E$-classes, it is injective, and therefore a
	bijection on structures of size $n$.
\end{proof}

\begin{lemma} \label{lemma:reduct}
Let $\mathcal{L}$ be a language and $\mathcal{L}'$ a reduct of $\mathcal{L}$.
Given a class $\mathcal{C}$ of $\mathcal{L}$-structures which admits a 
logical limit law, any class $\mathcal{C}' $ of $\mathcal{L}'$-structures 
which expand uniquely to $\mathcal{C}$-structures also admits a logical limit law. 
\end{lemma}

\begin{proof}
The map $f$ is taken to be the map sending a structure in $\mathcal{C}'$ to
its unique expansion in $\mathcal{C}$.  Because this expansion is unique, $f$ is
bijective on structures of size $n$ for all $n$.  We take $g$ to be the identity
map on formulas (as $\mathcal{L}'$ is a reduct of $\mathcal{L}$).  Then these maps satisfy the requirements of Lemma \ref{lemma:translation}.
\end{proof}

\begin{lemma}\label{lemma:expand}
Every composition expands uniquely to a fractured order, up to
isomorphism.
\end{lemma}

\begin{proof}
There is a unique way (up to isomorphism) to linearly order
each $\mathrel{E}$-class individually.  Because ordering these classes determines
$\prec_2$, there is a unique (up to isomorphism) way to define $\prec_2$ on any
composition, expanding it to a fractured order.
\end{proof}

\begin{theorem}\label{theorem:oeqlimlaw}
Compositions admit an unlabeled logical limit law.
\end{theorem}

\begin{proof}
By Lemma \ref{lemma:c-oer}, fractured orders are uniformly interdefinable
with convex linear orders (which admit an unlabeled logical limit law).  
Because every unlabeled composition expands uniquely to a 
fractured order, by Lemma \ref{lemma:reduct} we have an unlabeled limit law for
compositions.
\end{proof}

\section{Questions}
Our methods prove only an unlabeled limit law for compositions, hence, it is natural to ask the following:
\begin{question}
Do compositions admit a labeled logical limit law?
\end{question}

The original motivation for this work was to consider limit laws for permutation
classes. 
Compton \cite{Bell} devised a method for proving limit laws on
classes of structures based on analyzing growth rate, i.e., function
counting the number of unlabeled structures of size $n$ in the class, for each $n$, assuming the classes are
closed under disjoint union.  Although linearly ordered structures are not
closed under disjoint union, $\oplus$ provides a non-symmetric analogue.

\begin{question} \label{question:compton}
Can the method of Compton outlined in \cite{Bell} be extended to classes of
ordered structures using the $\oplus$ operation?
\end{question}

In particular, Compton's method shows that slow growth rate gives rise to
limit laws.  Although the class of all permutations has
growth rate $n!$, any other class has at most exponential growth rate \cite{marcus2004excluded}, which
is comparatively slow.

A positive answer to Question \ref{question:compton} would suggest the following, 
although it might be approached by other means.
\begin{question}
Do all sum-closed permutation classes (besides the class of all permutations) 
admit a logical limit law?
\end{question}


The gap in growth rates between proper permutation classes and the class of all permutations is one manifestation of the fact that proper permutation classes are comparatively tame. This is also witnessed by the fact that they have bounded twin-width, in the sense of \cite{bonnet2021twin}, which may help answer the following.

\begin{question}
Do all permutation classes (besides the class of all permutations) admit a 
logical limit law?
\end{question}

\subsection{Acknowledgements}
We thank the referees for their comments improving the exposition of the paper.

\bibliography{limlaws}
\bibliographystyle{plain}

\end{document}